\documentclass[11pt]{amsart}
\usepackage{amsmath, amssymb, amsthm,color}
\usepackage{amssymb,longtable}
\usepackage{amsfonts}
\usepackage{setspace}
\usepackage{amsmath} 
\usepackage{hyperref} 
\usepackage{latexsym}
\usepackage{array}
\usepackage{amssymb}
\usepackage{showkeys}
\theoremstyle{plain}

\numberwithin{equation}{section}
\newtheorem{thm}{Theorem}[section]
\newtheorem{prop}[thm]{Proposition}
\newtheorem{lemm}[thm]{Lemma}

\theoremstyle{remark}

\newcommand{\bq}{\begin{equation}}
\newcommand{\eq}{\end{equation}}
\begin{document}

\title[Slightly Supercritical 2-D Euler Equations]{Osgood's Lemma and Some Results for the Slightly Supercritical 2D Euler Equations for Incompressible Flow}

\author{Tarek Mohamed Elgindi}
\address{The Courant Institute,\\ New-York University\\
251 Mercer Street\\   USA}
\email{elgindi@cims.nyu.edu}

\thanks{T.M. Elgindi is partially supported by NSF grant DMS-0807347.}

\maketitle
\begin{abstract}
We investiage the (slightly) super-critical 2-D Euler equations. The paper consists of two parts. In the first part we prove well-posedness in $C^s$ spaces for all $s>0.$ We also give growth estimates for the $C^s$ norms of the vorticity for $0< s \leq 1.$  In the second part we prove global regularity for the vortex patch problem in the super-critical regime.This paper extends the results of Chae, Constantin, and Wu where they prove well-posedness for the so-called LogLog-Euler equation. We also extend the classical results of Chemin and Bertozzi-Constantin on the vortex patch problem to the slightly supercritical case. Both problems we study are done in the setting of the whole space.
 
\end{abstract}

\section{Introduction}

\setstretch{1.49}

The question of well-posedness for active scalar equations is one which has attracted much attention in the last several years. One of the breakthroughs in the study of active scalar equations was the proof(s) of well-posedness for the 2-D critically dissipative Surface Quasi-Geostrophic (SQG) equations.  A proof of wellposedness was discovered by Kiselev, Nazarov, and Volberg using their  \textquotedblleft modulus of continuity method\textquotedblright [14]. Around the same time, Caffarelli and Vasseur proved the well-posedness for SQG by extending De-Giorgi's method for non-linear elliptic equations to non-local equations [5]. A third proof was discovered by Kiselev and Nazarov using the duality of $H^p$ and $C^{\alpha}$ for $p<1$ [15]. A fourth proof was discovered by Constantin and Vicol using their \textquotedblleft Nonlinear Maximum Principle\textquotedblright method [9]. Finally, yet another proof has been put forth by Maekawa and Miura by applying Nash's work on non-linear elliptic equations to non-local equations [19]. All of this was done in the   \textquotedblleft critical  \textquotedblright regime. In the super-critical regime, we have some negative results and very few positive results. The classic negative result is the emergence of shocks in the inviscid Burgers equations. Kiselev, Nazarov, and Shterenberg showed that shocks may emerge for the viscous Burgers equation for dissipation below the critical level [16]. Another proof of this result is that of Dong, Du, and Li [13]. See [8], [10], and [18] for further blow-up results on a one-dimensional model of SQG.   
It is said that the first \emph{positive} result in the super-critical regime for fluid equations is that of Tao [20] where he proved well-posedness for the slightly super-critical hyperviscous 3D navier stokes equations. Chae, Constantin, and Wu used Littlewood-Paley theory to prove well-posedness for the supercritical LogLog-Euler [6]. Dabkowski et al. have recently published a work on a  \textquotedblleft slightly super-critical  \textquotedblright SQG model--where they use the modulus of continuity method to move slightly into the supercritical regime [12].  In this paper, we give a proof of well-posedness for a slightly more singular equation than LogLog-Euler.  
\\
The main question we address is: how far into the super-critical regime can we prove well-posedness for inviscid active scalar equations using the  \textquotedblleft conventional  \textquotedblright methods? We ask the same question for vortex patch problems as well. We prove here that so long as the velocity has an Osgood modulus of continuity, we have well-posedness. This is reminiscent of the uniqueness theorems of Yudovich in [22,23] and Vishik in [21] for the 2-D Euler equations. See also the result of Bahouri and Chemin [1]. In [14] Kelliher hints at the question: how large is the class of initial data for which uniqueness for 2-D Euler holds? Kelliher speculates that the answer to this question may be related to whether the modulus of continuity of the velocity is Osgood.   
\

We first recall a version of Osgood's lemma [2].

\begin{lemm}

Let $\gamma$ be a positive, nondecreasing, continuous function defined on $[0,\infty)$.

Suppose $$\int_{2}^{\infty}\frac{dr}{\gamma(r)}= \infty.$$ 

Suppose further that $$f(t) \leq f(0)+\int_{0}^{t}\gamma(f(r))dr$$

Let $$H(r)=\int_{1}^{r} \frac{dr}{\gamma(r)}$$

Then $$f(t) \leq H^{-1}(H(f(0))+t).$$

\end{lemm}

\emph{Proof:} 

\vspace{2mm}

Let $R(t)= f(0)+\int_{0}^{t} g(r)\gamma(f(r))dr.$ 
Then $$R'(t)= g(t)\gamma(f(t)).$$
Which implies that
$$R'(t) \leq g(t) \gamma(R(t)).$$
Dividing by $\gamma(R(t))$ and integrating gives the lemma.
\qed

The following result is a corollary of Osgood's lemma.

\begin{lemm}

Let $\gamma$ be a positive, nondecreasing, continuous function defined on $[0,\infty)$. Suppose that $\gamma(x+y)\leq C(\gamma(x)+\gamma(y)).$

Suppose $$\int_{2}^{\infty}\frac{dr}{\gamma(r)}= \infty.$$ 

Suppose further that $$f(t) \leq f(0)+\int_{0}^{t}\gamma(f(r))dr+ \gamma(\int_{0}^{t}f(r)dr)$$

Let $$H(r)=\int_{1}^{r} \frac{dr}{\gamma(r)}$$

Then $$f(t) \leq H^{-1}(H(f(0))+C(t^2+t)).$$

\end{lemm}

The proof of this lemma relies on using the condition that $\gamma(x+y)\leq C(\gamma(x)+\gamma(y))$ and the fact that $\gamma$ cannot grow faster than, say, quadratic to prove that   

$$ \gamma(\int_{0}^{t}f(r)dr)\leq C(t+1)\int_{0}^{t}\gamma(f(r))dr.$$

\section{Part I: Well-posedness in H\"older Spaces}

\begin{thm}

Suppose $m$ is a non-decreasing, positive, measurable function on $[0,\infty).$ Suppose that $m$ satisfies $m(2t) \leq Cm(t)$ for all $t\in [0,\infty)$ and for some constant $C>0.$  
Suppose further that  \begin{equation} \int_{2}^{\infty}\frac{1}{t\text{Log}(t)m(t)}dt = +\infty \end{equation} 
Then the following system is globally well-posed in $C^s \cap L^2$ for $0<s\leq 1.$ 
$$\omega_t + u \cdot \nabla \omega=0.$$
$$u=m(|D|)\nabla^{\perp}\Delta^{-1}\omega$$

Set $\Gamma(s)=m(s)(1+\text{Log}(s)).$  Then we define $H(s)=\int_{1/a}^{s} \frac{dr}{r\Gamma(r)}.$
Finally, set $f(t)=\text{Log}(|\omega(t)|_{C^s \cap L^2},$ $0<s\leq 1.$
Then $\omega$ satisfies the following growth estimate:  $$f(t) \leq H^{-1}(H(f(0))+Ctf(0)).$$

As a consequence, we have global well-posedness in $C_{c}^{\infty}.$

\qed
\end{thm}

\emph{Remarks:} 

(1)Our proof also gives global well-posedness in higher dimensions as well.  
 \

(2)Our proof also goes through for the case where $div(u) \neq 0.$
\

 Assumption (2.1) is precisely what is needed to guarantee that the velocity $u$ has (a-priori) an Osgood modulus of continuity.

Our focus here is on pointing to the heart of the matter, which is the Osgood lemma. 

The best result that existed prior to the writing of this is that of Chae, Constantin, and Wu [6]: wellposedness in Besov spaces in the case where $m(|\xi|)= \text{Log}(1+\text{Log}(|\xi|^2+1))^\gamma,$ where $\gamma \in [0,1]$. It is possible to show well-posedness in general Besov spaces using the method that we describe--we do not do this for the sake of keeping our arguments transperant. Of course, the result of Chae, Constantin, and Wu can be subsumed under the Besov space version of Theorem 1.1. Unfortunately Theorem 1.1 is not much better than the results of [6]. Indeed, examples 
of super-critical regimes for which our theorem gives well-posedness are:
\\
(A) $m(|\xi|)= \text{Log}(1+\text{Log}(|\xi|^2+1))^{\gamma},$ which is the same as what was given in [6].
\\
(B)$m(|\xi|)= \text{Log}(1+\text{Log}(|\xi|^2+1))^{\gamma_1}...\text{Log}(1+\text{Log}(1..+\text{Log}(|\xi|^2+1)))^{\gamma_n}$

Where $\gamma, \gamma_1,...\gamma_n \in [0,1].$

We remark that there is another proof of this result by Dabkowski et al. in [11].

We prove Theorem 1.1 by proving a certain logarithmic-Sobolev type inequality and then we rely upon the following Beale, Kato, and Majda-type result:

\

Let $\omega$ be a $C^{s}\cap L^2 $  solution of the system in Theorem 1.1 on $[0,T],$ with $0<s\leq 1,$

then, we have that
\begin{equation} |\omega(T)|_{C^s \cap L^2} \leq |\omega_0|_{C^s \cap L^2}\exp({\int_{0}^{T} |\nabla u(t)|_\infty dt}).\end{equation}
(2.2) follows from noting that $\omega(t)=\omega_{0}(\Phi(t))$ where $\Phi$ is the Lagrangian flow determined by $u$ from the system of ordinary differential equations: $$\dot{\Phi}=u(\Phi),$$ $$\Phi(0)=\text{Id}.$$  (2.2) is a consequence of the fact that $|\Phi|_{\text{Lip}}(T) \leq \exp(\int_{0}^{T}|\nabla u|_{L^\infty}(\tau)d\tau).$

Therefore the main step in proving global well-posedness is to bound $|\nabla u|_\infty .$ We will show that we can bound $|\nabla u|_\infty$ in terms of the $|\omega|_{C^s \cap L^2}$ in such a way that the estimate in (2.2) can be closed using Osgood's lemma.

\subsection{The Littlewood-Paley Decomposition}
\

A very convenient tool for studying the action of Fourier multipliers is the Littlewood-Paley decomposition. The decomposition allows us to restrict our attention to looking at functions whose Fourier transform has support in annuli. Such a study begins with the following basic proposition.

\begin{prop} 
\

There exists a pair of smooth radial functions $\phi$ and $\chi$ whose support is in $B(0,2)-B(0,\frac{1}{2})$ and $B(0,1)$ respectively which satisfy

$$ \chi(\xi) + \sum_{j=0}^{\infty} \phi(2^{-j}\xi)=1$$

$$ |j-k| >1 \rightarrow \phi(2^{-j}\cdot) \cap \phi(2^{-k}\cdot) = \emptyset$$
\qed
\end{prop}

We define the following operators: 

$$\Delta_{-1}u= \mathcal{F}^{-1}(\chi(\xi) \mathcal{F}u)$$
$$\Delta_{j}u= \mathcal{F}^{-1}(\phi(2^{-j}\xi) \mathcal{F}u)$$

We may decompose any tempered distribution $f \in \mathcal{S}'(\mathbb{R}^{n})$ as $$f=\sum_{j=-1}^{\infty} \Delta_{j}f.$$

Define the operator $S_{j}$ in the following way:

$$S_{j} =\sum_{k=-1}^{j}\Delta_{k}f.$$

That is, $S_j f$ is a cut-off at frequency $2^j$ of $f.$

It is easy to show that $$|S_{j} m(|D|)f|_{L^\infty} \leq m(2^j)|f|_\infty ,$$ if $f$  is a Schwartz function and $m$ satisfies mild conditions (such as those in Lemma 2.5). The proof of this can be seen by noting that if $\psi$ is a smooth function then $|m(|D|)\psi(\lambda \, \cdot)|_{L^\infty}\lesssim m(\lambda). $

Now define the following norms:

$$|u|_X \equiv ({\sum_{j=-1}^{\infty} 2^{2js}{|\Delta_j u|_{L^2}}^2})^{\frac{1}{2}}$$
and 
$$|u|_Y \equiv {\sup_{j}} \, 2^{js}{|\Delta_j u|_{L^\infty}}.$$

Define the following spaces:

$$X\equiv \{ f \in \mathcal{S}'(\mathbb{R}^{n}): |u|_{X}<\infty\}$$

$$Y\equiv \{ f \in \mathcal{S}'(\mathbb{R}^{n}): |u|_{Y}<\infty\}$$

Where $\mathcal{S}'$ denotes the space of all tempered distributions. We remark that the space $X$ defined here is the same as the Besov space $B^{s}_{2,2}$ and $Y$ is the Besov space $B^{s}_{\infty,\infty}$  For more on Besov spaces see [2].

We note here that the particular choice of $\phi$ and $\chi$ have no bearing on the function spaces we have defined. 
\

\begin{prop} 

Let $H^s$  denote the Sobolev space of index $s\geq 0$ on $\mathbb{R}^{n}.$ Then $$H^s = X$$ and there exists a constant $C=C_{s}>0$ so that \begin{equation}\frac{1}{C} |u|_{H^s} \leq |u|_X \leq C |u|_{H^s}.  \end{equation} Let $0<s<1.$ Then $$Y\equiv C^{s}$$ and there exists a universal constant $C=C_{s}>0$ so that  \begin{equation}\frac{1}{C} |u|_{C^s} \leq |u|_Y \leq C |u|_{C^s}.  \end{equation}
\qed
\end{prop}

\begin{prop}
\

(A) $|f|_{\infty} \leq  \sum_{j=-1}^{\infty} |\Delta_{j}f|_{\infty}.$
\

(B) Let $R$ be a Calderon-Zygmund operator. Then $|\Delta_j R(f)|_{p} \leq C|\Delta_j f|_{p}, 1\leq p \leq \infty,$ where $C$ is independent of $p$ and $j\geq 0$
\qed
\end{prop}

\subsection{The Main Inequality for Theorem 1.1}
\

\begin{lemm}
\

Let $m$ be a a positive, non-decreasing, measurable function on $[0, \infty).$ Suppose that $m(2x)\leq Cm(x)$ for some constant $C>0$ and all $\xi\in \mathbb{R}^{2}$. Suppose further that $|m(|\xi|)| \leq C \text{Log}(|\xi|+2).$
.Suppose that $f$ and $g$ are Schwartz functions. Suppose further that $\hat{f}(\xi)=m(|\xi|)R(\xi)\hat{g}(\xi)$ where $R$ is any Calderon-Zygmund operator (such as a Riesz transform or a composition of Riesz transforms). Let $s>0.$
Then the following inequality holds:

\begin{equation} |f|_{\infty} \leq C|g|_{L^2} + C_s |g|_{\infty}(1+\text{Log}(\frac{|g|_{C^s}}{|g|_{_\infty}})m(\frac{|g|_{C^s}}{|g|_{_\infty}})). \end{equation}
\qed
\end{lemm}
\emph{Proof:}
Let $s'<s.$
$$\sum_{-1}^{\infty} |\Delta_j f |_\infty = |\Delta_{-1} f|_{\infty} +\sum_{0}^{N} |\Delta_j f |_\infty+ \sum_{N+1}^{\infty} |\Delta_j f |_\infty$$
$$ \leq C|g|_{L^2}+C \cdot N\cdot m(2^{N+1}) \sup_j|\Delta_j g|+\sum_{N+1}^{\infty} |\Delta_j f |_\infty$$
$$ \leq C|g|_{L^2} + C \cdot N\cdot m(2^{N+1}) |g|_{\infty} + C \sum_{N+1}^{\infty} 2^{-js'}2^{js'} m(2^j)|\Delta_j g |_\infty$$
We now use the cauchy-schwarz inequality on the second term and get:
$$ \leq C|g|_{L^2} + C \cdot N\cdot m(2^{N+1}) |g|_{\infty} + C_{s}2^{(-N-1)s'} \sum_{N+1}^{\infty}2^{sj}  |\Delta_j g |_{\infty}$$
Now by Proposition 2.4 we get:
$$=C|g|_{L^2} + C \cdot N\cdot m(2^{N+1}) |g|_{\infty} + C(s)2^{(-N-1)s'}|g|_{C^s}.  $$

Taking $N \approx \text{Log} (\frac{|g|_{C^s}}{|g|_{\infty}})$ finishes the proof.
\qed
\\
\\

\subsection{Proof of Theorem 1.1}

Let $s>1.$ By (2.2), we have:

$$|\omega(t)|_{C^s\cap L^2} \leq |\omega_{0}|_{C^s\cap L^2} \exp({\int_{0}^{t} |\nabla u(\tau)|_{\infty}d\tau})$$ 
\
\
Now use (2.4) with $f=\nabla u$ and $g=\omega:$

$$ |\omega(t)|_{C^s\cap L^2} \leq |\omega_{0}|_{C^s\cap L^2} e^{C \int_{0}^{t}|\omega_{0}|_{L^2} + |\omega_{0}|_{\infty}(1+\text{Log}(\frac{|\omega(\tau)|_{C^s\cap L^2}}{|\omega_{0}|_{\infty}})m(\frac{|\omega(\tau)|_{C^s\cap L^2}}{|\omega_{0}|_{\infty}}))d{\tau}},$$

\
$$ \text{Log}(|\omega(t)|_{C^s\cap L^2}) \leq \text{Log}(|\omega_{0}|_{C^s\cap L^2} )+{C \int_{0}^{t}|\omega_{0}|_{L^2}+|\omega_{0}|_{\infty}(1+\text{Log}(\frac{|\omega(\tau)|_{C^s\cap L^2}}{|\omega_{0}|_{\infty}})m(\frac{|\omega(\tau)|_{C^s\cap L^2}}{|\omega_{0}|_{\infty}}))d{\tau}}.$$
\

Note that we have used the fact that $|\omega_{0}|_{L^p}$ is conserved for $ 1<p\leq \infty .$ We now recall that $|\omega_0|_{\infty} \leq C|\omega_0|_{C^s\cap L^2}$ because $s>1.$ Then use the assumption (2.1) and Osgood's lemma to conclude that $$|\omega(t)|_{C^s\cap L^2} \leq M(t,|\omega(0)|_{C^s \cap L^2})$$
for some $M\in L^{\infty}_{loc}.$ This estimate along with a local well-posedness result gives that any local solution may be continued indefinitely which gives the global well-posedness.
 \qed

\subsection{Growth Estimates}

It is known that the best growth estimates for $C^s$ norms of $\omega$ for the 2-D Euler equations are of double exponential type. 
We can generalize this result to our slightly supercritical equations. Indeed, using Osgood's lemma as is written in Section 1 we first set $\Gamma(r)=m(r)(1+\text{Log}(r)).$ Then we define $H(r)=\int_{1/a}^{r} \frac{dr'}{r'\Gamma(r')}.$
Finally, set $f(t)=Log(|\omega(t)|_{C^s}), \, 0<s\leq 1.$ 
Then we $\omega$ satisfies the following growth estimate:  $$f(t) \leq H^{-1}(H(f(0))+Ctf(0))$$ using the estimates above.
We easily see that if, for example, $$m(|\xi|)= \text{Log}(1+\text{Log}(|\xi|^2+1))^{\gamma_1}...\text{Log}(1+\text{Log}(1..+\text{Log}(|\xi|^2+1)))^{\gamma_n},$$ $\gamma_{i}=1, \ i=1...n,$ then we get that $H^{-1}$ is like a composition of $n+2$ exponentials to that $|\omega|_{C^s}$ is bounded by $$\exp(\exp(....(\exp(Ct|\omega_0|_{C^s\cap L^2})))),$$ where $n+2$ exponentials are taken. Proving that such growth estimates are attainable would, of course, reveal quite a bit about the supercritical regime.

\subsection{Local Well-Posedness}

Local well-posedness may be shown by considering a "splitting" approximation to the Euler system:

$$\omega^{(0)}=\omega_0 ,$$

$$u^{(n)}=m(|D|)\nabla^{\perp}\Delta^{-1}\omega^{(n-1)}, \, n\geq 1,$$

$$\omega^{(n)}_{t} + u^{(n)} \cdot \nabla \omega^{(n)}=0, \, n\geq 1.$$

Using the a-priori estimates of the previous sections we can easily get local well-posedness and then global well-posedness will follow. 
 
\subsection{Higher Regularity} The above estimates give global well-posedness for $\omega$ in $C^s$ for $0<s\leq 1.$ Higher regularity (say for data in $C_{c}^{\infty}$) follows easily by standard estimates since the Lipschitz norm of both $u$ and $\omega$ are under control. 

\subsection{Further Applications} 

We remark finally that Lemma 2.5 has more applications than what we have presented above. In particular, Lemma 2.5 can be used to derive sharper continuation criteria for the 3-D Euler equations. This, along with other applications, will be elucidated in a forthcoming paper. 

\section{Part II: The Vortex Patch Problem}

This section is devoted to proving global-wellposedness for what we will call slightly singular vortex patch problems. Recall the vortex patch problem as was studied in [4] and [7]:

$$\phi_t +u\cdot \nabla \phi=0,$$
$$\phi(x,0)=\phi_{0}(x),$$
$$v(x,t)= \frac{\omega_{0}}{2\pi}\int_{E}\nabla^{\perp}_{x} \text{Log}|x-y|dy.$$

Where $$E = \{ x\in \mathbb{R}^2 | \, \phi(x) >0 \}.$$ 
Assume

$$E_{0} = \{ x\in \mathbb{R}^2 | \, \phi_{0}(x) >0 \}$$ is bounded and has a smooth boundary. In particular, we assume that $$ \inf_{x\in \partial E_{0}} |\nabla \phi_{0}(x)| \geq l>0$$ and $\phi \in C^{1,\mu}$ for some $\mu \in (0,1].$

Chemin showed global regularity of solutions to the vortex patch problem with smooth initial data [7]. Bertozzi and Constantin later simplified Chemin's argument in [4]. Here, we combine the approaches of Bertozzi-Constantin and Chemin to include a more singular expression for $v$ than what is given by the Biot-Savart law. The main steps in both approaches are more or less the same. Both approaches aimed to show the propogation of the $C^{1,\mu}$ regularity of the boundary. Towards this aim, the authors first got a bound on the Lipschitz norm of the velocity field in terms of the boundary regularity of the vortex patch (i.e. the regularity of $\phi$). We intend to study a more singular problem than the one studied in [4] and [7]. Due to this extra singularity, the Lipschitz norm of the velocity field will not be bounded which will introduce several extra difficulties. To overcome these difficulties we will need some new ideas.  

Indeed, we show global regularity for the following system:

$$\phi_t +u\cdot \nabla \phi=0,$$
$$\phi(x,0)=\phi_{0}(x),$$
$$u(x,t)= \frac{a_{0}}{2\pi}m(|D|)\int_{E}\nabla^{\perp}_{x} (\text{Log}|x-y|)dy.$$ 
We call the above formula for $u$ the Modified Biot-Savart Law.

\begin{thm}

Suppose $m$ is a non-decreasing, positive, measurable function on $[0,\infty).$ Suppose that $m$ satisfies $m(|\xi_1||\xi_2|) \leq C(m(|\xi_1|)+m(|\xi_2|))$ for all $\xi_1, \xi_2 \in \mathbb{R}^{2}$ and for some constant $C>0.$ Also suppose that $m(|\xi|)\leq C { \text{Log}(|\xi|+2)}.$  Finally (and most importantly), we assume $m$ satisfies the following:

$$\int_{2}^{\infty} \frac{dr'}{r'm(r')(1+Log(r'))}=+\infty.$$

Given $\omega_0 \neq 0,$ $E_0$ bounded and $\phi_0 \in C^{1,\mu}$ satisfying $ \inf_{x\in \partial E_{0}} |\nabla \phi_{0}(x)| \geq l>0,$ the solution to the modified vortex patch problem exists for all time. Moreover, there exists $C$ depending upon $|\omega_0|, area(E_0), |\nabla \phi_0|_{\mu}, |\nabla \phi_0|_{\infty},$ and $|\nabla \phi_{0}|_{\inf}$ so that  given $\epsilon>0$ the following estimates hold:

$$|\nabla v(t)|_{\infty} \leq   |\nabla v(0)|_{\infty}H^{-1}(tC({\epsilon}))$$
$$|\nabla \phi(t)|_{\mu-\epsilon} \leq   |\nabla \phi(0)|_{\mu}\exp(H^{-1}(tC({\epsilon})))$$
$$|\nabla \phi(t)|_{\infty} \leq   |\nabla \phi(0)|_{\infty}\exp(H^{-1}(tC({\epsilon})))$$
$$|\nabla \phi(t)|_{\inf} \geq   |\nabla \phi(0)|_{\infty}\exp(-H^{-1}(tC({\epsilon})))$$

Where $H$ is defined in terms of $m$ as in Section (4): $$H(r)=\int_{2}^{r} \frac{dr'}{r'm(r')(1+\text{Log}(r'))}.$$
\qed
\end{thm}

\begin{thm}

Suppose $m$ is a non-decreasing, positive, measurable function on $[0,\infty).$ Suppose that $m$ satisfies $m(|\xi_1||\xi_2|) \leq C(m(|\xi_1|)+m(|\xi_2|))$ for all $\xi_1, \xi_2 \in \mathbb{R}^{2}$ and for some constant $C>0.$ Also suppose that $m(|\xi|)\leq C \text{Log}(|\xi|+2).$  Finally (and most importantly), we assume $m$ satisfies the following:

$$\int_{2}^{\infty} \frac{dr'}{r'm(r')(1+\text{Log}(r'))}=+\infty.$$  
Given $\omega_0 \neq 0,$ $E_0$ bounded and $\phi_0 \in C^{\infty}$ satisfying $ \inf_{x\in \partial E_{0}} |\nabla \phi_{0}(x)| \geq l>0,$ the solution to the modified vortex patch problem exists for all time. Moreover $\phi(t) \in C^{\infty}.$

\qed
\end{thm}

We assume, for simplicity, that $area(E_0)=area(E)=1.$

We are going to break the proof into estimating two quantities which will control the $C^{1,\mu}$ regularity of the boundary:

\vspace{2mm}

(1) $|\nabla \phi|_{\inf},$

\vspace{2mm}

(2)$|\nabla \phi|_\mu .$

\vspace{2mm}

\subsection{Well-posedness for the Classical Vortex Patch Problem}

The classical method, as given in the work of Bertozzi-Constantin [4] are based on the following estimates:

\vspace{2mm}

(A) $|\nabla u|_{\infty} \leq C \text{Log}(1+\Delta_{\mu})$

\vspace{2mm}

(B) $|\nabla u \nabla^{\perp} \phi|_{\sigma} \leq C|\nabla u|_{\infty}|\nabla^\perp \phi|_{\sigma}.$

\vspace{2mm}

Proving (A) involves studying the Biot-Savart law. Indeed, in the classical case, $\nabla u$ is derived from $\omega$ by a singular integral operator. Therefore, since $\omega$ is only $L^\infty,$ it isn't immediately clear that $\nabla u$ will also be bounded. However, due to the fact that $\omega$ is the characteristic function of a set with a $C^{1,\mu}$ boundary, it turns out that $\nabla u$ can be shown to be bounded. In our case, not even $\omega$ is bounded--which is a great source of difficulty in the slightly super-critical case.  

Proving (B), on the other hand, uses the observation that $\nabla u$ is much more regular in the tangential direction than in the normal direction. Indeed, it is possible to show in the classical case that if $W$ is a divergence free vector field which is tangent to $\partial E,$ then $$\nabla u W(x)= \int_{E} \frac{\sigma(x-y)}{|x-y|^2}(W(x)-W(y)),dy$$ where $\sigma$ is a function which is homogeneous of degree zero and which has mean zero on half circles. 

This crucial observation, observed first by Chemin [7], is what allows us to propagate $C^{1,\mu}$ regularity of $\phi$ without $\nabla u$ being itself regular.

\par

Now with (A) and (B) in hand and noting that  $$\nabla^\perp \phi_t +(u\cdot \nabla) \nabla^\perp \phi= \nabla u \nabla^\perp \phi,$$ we can prove the following two estimates:  

\begin{equation} |\nabla^\perp \phi|_{\inf} \geq \exp(-\int_{0}^{t}|\nabla u|_{\infty}), \end{equation} 
and
\begin{equation} |\nabla^\perp \phi|_{C^\mu} \leq |\nabla^\perp \phi_0 |_{C^\mu}\exp(\int^{t}_{0} |\nabla u|_\infty). \end{equation}

\subsection{The Difficulties With the Extra Singularity}

As in the case of global well-posedness for $C^s$ solutions, adding the extra loglog singularity does not change the general strategy of the proof. However, some new methods are needed to deal with the extra singularity for several reasons. The main difficulty is that $\nabla u$ is no longer bounded. There is a general principle for transport-type equations: if the velocity field is not Lipschitz then it is not possible to propagate regularity. And, indeed, we are unable to show propagation of $C^{1,\mu}$ regularity. However, with the aid of the method of losing estimates, we can show that if at $t=0,$ $\phi$ is $C^{1,\mu}$ regular, then it remains $C^{1,\mu-\epsilon}$ regular for all time. Another issue we have to face is that we can no longer bound $|\nabla \phi|_{\inf}$ so easily as before (estimate 3.1). In order to overcome the fact that $\nabla u$ is not bounded, we discover that in order to bound $|\nabla \phi|_{\inf}$ it suffices to bound $\nabla u$ near $\partial E$ in the direction tangent to $\partial E.$ 
  
\section{Main Ideas that go into the Proof of Theorem 3.1}

For active scalar equations the main thing one wants to do is bound the Lipschitz norm of the velocity field. In the first problem we studied above, we were able to achieve such a bound using a generalized logarithmic-Sobolev inequality. For our vortex patch problem, we are unable to get such a bound on the velocity field. Indeed, because of the extra singularity (given by $m$), we are unable to say that the anti-symmetric part of $\nabla u$ is bounded ( the curl$(u)$ unbounded). Indeed, such a bound shouldn't even be expected. Therefore, we are only able to bound $\nabla u$ in a space just below $L^\infty.$ Due to this, we are unable to propogate the $C^{1,\mu}$ regularity of the vortex patch. Due to our inability to bound $\nabla u$ in $L^\infty ,$ we face two issues:

\vspace{2mm}

(1) \, Without a Lipschitz velocity field, the $C^{\alpha}$ regularity cannot be propogated.

\vspace{2mm}

(2) \, The bound from below on $|\nabla \phi|_{\inf}$ is lost.

\vspace{2mm}

The first issue (1), luckily, can be resolved by immitating the losing estimates of Bahouri and Chemin [2]. Indeed, if the velocity field is not Lipschitz we are unable to propagate $C^\alpha$ regularity, but if the velocity field is not too bad (say Osgood continuous), then we will be able to say that the regularity of $\phi$ will remain $C^{\alpha-\epsilon}$ if we start with $C^\alpha$ regularity. 
The second issue (2), unfortunately, cannot be circumvented by anything like the losing estimates. Fortunately, we are able to recast our problem in such a way that even though the antisymmetric part of $\nabla u$ is unbounded,  $\nabla u$ actually remains bounded near $\partial E$ in the direction tangent to $\partial E$. 

To get the results of our theorem will rely upon a series of estimates. 

\vspace{5mm}

\emph{Step 1:}

\vspace{5mm}

Consider $\nabla u.$ Let $$\Delta_\mu = \frac{|\nabla \phi|_\mu}{|\nabla \phi|_{\inf}}$$

The following estimate holds:
\begin{equation}  |S_j\nabla u|_{\infty} \leq C(\mu)m(2^j)(1+\text{Log}(\Delta_\mu)). \end{equation}

We get this estimate by writing the symmetric part of $\nabla u,$ $S(u),$ in the following way:

$$S(u)(x)= m(|D|)\int_{E} \frac{\sigma(x-y)}{|x-y|^2}dy$$ where $$ \sigma(z)=\frac{1}{|z|^2}\left( \begin{array}{cc}
-2(z_2 )(z_1 ) & (z_1)^2-(z_2 )^2 \\
(z_1 )^2 -(z_2 )^2 & 2(z_2 )(z_1 ) \\
\end{array} \right).$$ 

Furthermore, we will show the follwing estimate which says that the tangential part of $\nabla u$ is bounded even though $\nabla u$ itself is not:

\begin{equation} |\nabla u \nabla^\perp \phi \cdot \nabla^\perp \phi(x)  | \leq Cm(\Delta_\mu) (1+\text{Log}(\Delta_\mu))|\nabla^\perp \phi(x)|^2 . \end{equation}

Of course this estimate was obvious in the classical case when $\nabla u$ was bounded. However, since $\nabla u$ is no longer bounded, we are not able to get this estimate so easily.

\vspace{5mm}

\emph{Step 2:}

\vspace{5mm}

Apply $\nabla^{\perp}$ to the transport equation:

$$ \phi_t +u\cdot \nabla \phi=0$$

\begin{equation} \nabla^{\perp}\phi_t +(u\cdot \nabla) \nabla^{\perp} \phi= \nabla u \nabla^{\perp}\phi \end{equation}

We want to get a good estimate of the non-linearity on the right-hand side of (4.2).  Notice that $\nabla^{\perp}\phi$ is tangent to the boundary of $E$ and is divergence free. 

\vspace{5mm}

\emph{Step 3:}

\vspace{5mm}

If $u$ is given by the modified Biot-Savart law and if $W$ is a divergence free vector tangent to $\partial E$, then the following holds:

\begin{equation} \nabla u W(x)= m(|D|)\int_{E}\frac{\sigma(x-y)}{|x-y|^2}(W(x)-W(y))dy + C_1 \end{equation} 
\begin{equation} \nabla u W(x)= \int_{E}\sigma(x-y)G(|x-y|)(W(x)-W(y))dy + C_1 \end{equation} where $G(|x-y|) \leq \frac{C}{|x-y|^2} (1+m(\frac{1}{|x-y|})).$

$C_1$ is such that the estimate in Step 4 is valid. 

\vspace{5mm}

\emph{Step 4:}

\begin{equation} \sup_{j}2^{j\mu}\frac{1}{m(2^j)}|\Delta_j (\nabla u W)|_{\infty} \leq C(|S(u)|_\infty +1)|W|_{\mu},  \end{equation} for any $\mu>0.$

\emph{Step 5:}

\vspace{5mm}

Apply $\Delta_j$ to (4.2). Then we get:

$$ \Delta_j (\nabla^{\perp} \phi )_t +(u\cdot \nabla) \Delta_j (\nabla^{\perp} \phi )= \Delta_j(\nabla u \nabla^{\perp} \phi) +CO_j$$

Where $CO_j$ is the following commutator:

$$CO_j = [u\cdot \nabla, \Delta_j] \nabla^{\perp}\phi.$$

Using the bounds from Step 4 above, $$  2^{\sigma j}|CO_j|_{\infty} \leq  2^{\sigma j}(C(\mu)(1+m(\Delta_{\mu})\text{Log}(\Delta{\mu})) |\Delta_j (\nabla^{\perp}\phi)|_{\infty},$$ with constants independent of $j.$

\vspace{5mm}

\emph{Step 6:}

\vspace{5mm}

We need a good estimate on $|\nabla \phi|_{\inf}.$
Take (4.2) and dot it with $\nabla^\perp \phi.$ Now divide by $|\phi|^2.$ Then we get:

$$(\text{Log}|\nabla \phi|)_t = \nabla u \nabla^\perp \phi \cdot \frac{\nabla^\perp \phi}{|\nabla^\perp \phi|^2}.$$  

Note now that the right hand side above is bounded because the tangential part of $\nabla u$ is bounded.

We then see that $$ |\nabla \phi|_{\inf} \geq |\nabla \phi_{0}|_{\inf} \exp(\int_{0}^{t} \text{Log}(1+\Delta_\mu)m(\Delta_\mu)(\tau)d\tau).$$ 

\emph{Step 7:}

The main result of this step is:

$$|\nabla \phi|_{\mu-\epsilon}(T) \leq 2|\nabla \phi_{0}|_\mu \text{exp}[\frac{C(\gamma)}{\epsilon}m(\text{exp(V(T))})V(T)],$$

where, $V(t)=\int_{0}^{t}C(\gamma)(1+\text{Log}(\Delta_\gamma))(s)ds.$

This will follow from using the losing-estimate ideas of Bahouri-Chemin [1] to compensate for the $m(2^j)$ term in the bound on $\nabla u$ in Step 5, (4). We trade a little loss of regularity for the fact that $u$ isn't quite Lipschitz.

Once we have this bound, we may use Lemma 1.2 and the results of Steps 6 and 7 to conclude that the $C^{\mu-\epsilon}$ norm of $|\nabla \phi|$ is bounded for all finte time intervals for any $\epsilon>0.$

\section{Some Basic Tools}

In all that follows, we will distinguish between the actual velocity field $u$ and the "classical" velocity field $\tilde{u}:$ 

$$ u(x):= m(|D|) \int_{E} \nabla^\perp \text{Log}|x-y|dy$$
$$\tilde{u}(x):= \int_{E} \nabla^\perp \text{Log}|x-y|dy$$

First we want to understand the action of $m$ on the Fourier side as a convolution in physical space.

\begin{lemm}

$$S(u)(x)= \int_{E} \sigma(x-y)G(|x-y)|dy$$ where $$|G(|x-y|)| \leq C(1+m(\frac{1}{|x-y|}))\frac{1}{|x-y|^2}$$ and $$ \sigma(z)=\frac{1}{|z|^2}\left( \begin{array}{cc}
-2(z_2 )(z_1 ) & (z_1)^2-(z_2 )^2 \\
(z_1 )^2 -(z_2 )^2 & 2(z_2 )(z_1 ) \\
\end{array} \right).$$

\qed 
\end{lemm}

\vspace{5mm}

\emph{Proof of Lemma 5.1}

\vspace{5mm}

Consider the multiplier $$f(|\xi|)= \frac{m(|\xi|)}{|\xi|^2}.$$ Our aim is to compute $\hat{f}.$ Notice that because $f$ is radial, $\hat{f}$ will also be radial. We will write $\hat{f}(x)=\hat{f}(|x|)=\hat{f}(\rho)$ Recall also that since we are in two-dimensions, we have that $$\hat{f}(\rho)=\frac{1}{2\pi}\int_{0}^{\infty}J_{0}(2\pi \rho r)\frac{m(r)}{r}dr,$$ where $J_0$ is a Bessel function of the first kind defined by:

$$J_{0}(x)=\frac{1}{2\pi}\int_{-\pi}^{\pi}e^{-ixsin(\theta)}d\theta.$$

  Of course we encounter a problem with this integral at $r=0.$ To circumvent this problem we define the distribution $\frac{1}{|\cdot|^2}$ more properly as: $$<\frac{m(\cdot)}{|\cdot|^2}, g>\equiv \int_{B_1} \frac{(g(\xi)-g(0))m(|\xi|)}{|\xi|^2}d\xi+\int_{B_{1}^c} \frac{g(\xi)m(|\xi|)}{|\xi|^2}d\xi.$$ From this we get:
$$\hat{f}(\rho)=\frac{1}{2\pi}(\int_{0}^{1}(J_{0}(2\pi \rho r)-1)\frac{m(r)}{r}dr+\int_{1}^{\infty}J_{0}(2\pi \rho r)\frac{m(r)}{r}dr).$$
We now differentiate in $\rho$ $$\hat{f}'(\rho)=\int_{0}^{\infty}J_{0}'(2\pi \rho r)m(r) dr.$$ Notice here, by the way, that if $m \equiv 1$ this will imply that $\hat{f}(\rho)=C \text{Log}(\rho),$ which is the fundamental solution to Laplace's equation in two dimensions as is expected.
Taking one more derivative in $\rho,$ we get:
$$\hat{f}''(\rho)=2\pi\int_{0}^{\infty}J_{0}''(2\pi\rho r)rm(r)dr.$$ We note that the differentiations that we are taking are justified by integration by parts along with the $\frac{1}{\sqrt{x}}$ decay of $J_0$.
Now, by substitution $$\hat{f}''(\rho)=\frac{2\pi}{\rho^2} \int_{0}^{\infty}J_{0}''(2\pi r) r m(\frac{r}{\rho})dr.$$ Using that $m(xy) \leq C (m(x)+m(y)),$ we get $$\hat{f}''(\rho) \leq \frac{C}{\rho^2}(1+m(\frac{1}{\rho})).$$ This along with simple calculus concludes the proof of the lemma.
 
\qed

We will need the following standard commutator estimate.

\begin{lemm}

Let $f$ and $g$ be Schwartz functions. Let $m$ be a smooth function such that $m(|\xi|)\lesssim |\xi|,$ for $ |\xi| \geq 1.$ Then the following estimate holds:
$$\sup_j 2^{j\mu} \frac{1}{m(2^j)} |\Delta_j((m(|D|)f)g-m(|D|)(fg))|\leq C\sup_j 2^{j\mu}|\nabla f|_\infty |\Delta_j g|_\infty .$$

\end{lemm}

\section{The bound on the tangential part of $\nabla u$}

We want to prove that $\nabla u$ is bounded in a certain sense. Of course this is not true (even in the case where the boundary of $E$ is flat). At the boundary of $E,$ however, the vorticity isn't changing in the tangential direction (because outside of $E,$ the vorticity is $1$ and inside it is $0$). Therefore, we expect $u$ to be more regular in tangential directions near the boundary. Indeed, in this section, we prove that the tangential derivative of the tangential part of $u$ is bounded. More precisely, we take the tangent vectors to the boundary and extend them in a natural way to the whole space. Call the extension $\tau.$  We will show that $$<\nabla u \tau,\tau>$$ is uniformly bounded in space.  

\begin{lemm}
Assume $u$ is given by the Modified Biot-Savart law:

$$u(x,t)= \frac{a_{0}}{2\pi}m(|D|)\int_{E}\nabla^{\perp}_{x} (\text{Log}|x-y|)dy.$$

Where $$E = \{ x\in \mathbb{R}^2 | \, \phi(x) >0 \}.$$
Let $\tau=\frac{\nabla^\perp \phi}{|\nabla^\perp \phi|}.$ 
Then, $$(\nabla u \,  \tau) \cdot \tau \leq C(\mu)(1+m(\Delta_{\mu})\text{Log}(\Delta_\mu)).$$
\qed
\end{lemm}

Note that the quantity $\tau$ above is just the tangent to the boundary when $x_0$ approaches the boundary of $E$. This basically says that the tangential component of the gradient is bounded. Indeed, we will show later that this quantity actually has Holder regularity even though $\nabla u$ is unbounded. 

\subsection{Some Preliminary Computations}

We now investigate the matrix: $$ \sigma(z)=\frac{1}{|z|^2}\left( \begin{array}{cc}
-2(z_2 )(z_1 ) & (z_1)^2-(z_2 )^2 \\
(z_1 )^2 -(z_2 )^2 & 2(z_2 )(z_1 ) \\
\end{array} \right)$$ 

Notice that $\sigma$ satisfies the following properties:

\vspace{2mm}

(A) $\sigma$ is a smooth function and homogeneous of degree 0.
 
\vspace{2mm}

(B) $\sigma$ has mean zero on the unit circle

\vspace{2mm}

(C) $\sigma(z)=\sigma(-z).$

Remark that (A), (B) and (C) imply that $\sigma$ has mean zero on circles and half-circles. 
The geometric lemma of Constantin and Bertozzi tells us that if $E$ is sufficiently regular, then the intersection of $\partial E$ with a small circle looks like a half-circle. We do not use this here, but we felt it important to mention. 

\vspace{5mm}

\emph{Proof of Lemma 6.1}

\vspace{5mm}

Set $\Delta_{\mu}=\frac{|\nabla \phi|_{\mu}}{|\nabla \phi|_{\inf}}.$ Let $\delta$ satisfy $\delta^\mu=\frac{1}{2\Delta_\mu}.$ A consequence of this choice of $\delta$ is that if $x$ is within $\delta$ of the boundary, then $|\nabla\phi|(x)\geq \frac{1}{2}|\nabla \phi|_{\inf}$ (recall that in the definition of the quantity $|\cdot|_{\inf}$ , the infimum is taken only on the boundary). For $x_{0} \in E,$ define $d(x_0)$ as the distance from $x_0$ to $\partial E.$ Note, finally, that $<\nabla u\tau, \tau>=<S(u)\tau,\tau>.$

Recall that using Lemma 5.1, $S(u),$ the symmetric part of $\nabla u,$ can be written in the following way: 
 $$ S (u) (x_{0})= \frac{\omega_0}{2\pi} \int_{E \cap \{ |x_0 -y| \geq \delta \}} G(|x_0 -y|) \sigma(x_0 -y)dy$$ 
$$+\frac{\omega_0}{2\pi} \int_{E \cap \{ |x_0 -y| \leq \delta \}} G(|x_0 -y|) \sigma(x_0 -y) dy $$

$$= I_1 +I_2.$$

Changing to polar coordinates, it is clear that $I_1$ is bounded by $C(1+\text{Log}(\Delta_\mu)m(\Delta_\mu)).$  This is because $|\sigma| \leq 1$ and $m$ is a decreasing function. 

Now, if $d(x_0)>\delta,$ $I_2$ must vanish because $\sigma$ has mean zero on any circle. Therefore, we consider the case when  $d(x_0)<\delta.$ 

$$I_{2}=\frac{\omega_0}{2\pi} \int_{E \cap \{ |x_0 -y| \leq \delta \}} G(|x_0-y|) \sigma(x_0 -y) dy. $$

Now, without a loss of generality we may assume $x_0 =0.$ We may also assume that the closest point on the boundary to the origin is at $(0,d(x_0))$ so that the outward normal at that boundary point is (1,0). 

Now, since $\sigma$ has mean zero on circles (and half-circles), we see that we may remove a circle of size $d(x_0)$ around the origin from $I_2,$ so that $$I_{2}=\frac{\omega_0}{2\pi} \int_{E \cap \{ d(x_0) \leq |y| \leq \delta \}} G(|y|) \sigma(y) dy. $$ 

Now, we really want to estimate $(I_2 \tau) \cdot \tau.$ We now want to replace $t$ by the tangential vector at the boundary $(0,1).$ This will introduce an error. But we will be able to control it:

$$\tau-(1,0)=\frac{\nabla^\perp \phi}{|\nabla^\perp \phi|}(0,0)-\frac{\nabla^\perp \phi}{|\nabla^\perp \phi|}(0,d(x_0)).$$

Since $\phi$ is assumed to be $C^{1,\mu}$ and $|\nabla \phi|_{\inf}>0,$ we see that this quantity can be bounded in the following way: $$|\tau-(1,0)|=|\frac{\nabla^\perp \phi}{|\nabla^\perp \phi|}(0,0)-\frac{\nabla^\perp \phi}{|\nabla^\perp \phi|}(0,d(x_0))| \leq \frac{m(\Delta_\mu)}{m(\frac{1}{d(x_0)})} \, .$$

Indeed, let $f$ be a $C^{\mu}$ function. Suppose that $d(x_0) \leq \delta=|f|_{C^\mu}^{-\frac{1}{\mu}}.$

Then, $$|f(0,d(x_0))-f(0,0)|\leq d(x_0)^\mu |f|_{C^\mu}\leq \frac{1}{m(\frac{1}{d(x_0)})}\delta^\mu m(\frac{1}{\delta})|f|_{C^\mu}$$ $$=\frac{1}{m(\frac{1}{d(x_0)})}m(|f|_{C^\mu}).$$

The desired inequality on $\tau-(1,0)$ is then apparent once we note that $|\frac{\nabla^\perp \phi}{|\nabla^\perp \phi|}|_{\mu} \approx \Delta_{\mu},$ since $|\nabla \phi|\geq \frac{1}{2}|\nabla\phi|_{\inf}$ since $d(x_0)\leq \delta.$

The above computations tell us that we may replace $\tau$ by $(0,1)$ and the resulting \emph{error} which is smaller than:

$$|I_{2} \frac{m(\Delta_\mu)}{m(\frac{1}{d(x_0)})}|.$$

In addition we must control:

$$<I_2 \begin{pmatrix}1\\0\end{pmatrix} , \begin{pmatrix}1\\0\end{pmatrix}>
$$
 
We now want to break $I_2$ into two parts:

$$ \frac{\omega_0}{2\pi} \int_{P \cap \{ d(x_0)\leq |y| \leq \delta \}} G(|y|) \sigma(y) dy-\frac{\omega_0}{2\pi} \int_{E\Delta P \cap \{ d(x_0)\leq |y| \leq \delta \}} G(|y|) \sigma(y) dy $$
$$=I_{21}-I_{22}$$ 

Where $P$ is the region $\{ y_{2} \leq d(x_0) \}$ and $E\Delta P$ is the symmetric difference of $E$ and $P.$

\emph{Claim: }$$<I_{21} \begin{pmatrix}1\\0\end{pmatrix} , \begin{pmatrix}1\\0\end{pmatrix}>=0
$$

\emph{Proof of Claim:}

$$<I_{21} \begin{pmatrix}1\\0\end{pmatrix} , \begin{pmatrix}1\\0\end{pmatrix}>= \frac{\omega_0}{2\pi} \int_{P \cap \{ d(x_0)\leq |y| \leq \delta \}} \frac{y_{1}y_{2}}{|y|^2} G(|y|) dy
$$

Note that the integrand is odd in $y_{1}.$ Note further that the domain of integration is symmetric with respect to $y_{1}=0.$ Thus the integral is zero.

\qed

Now it remains to control $$<I_{22} \begin{pmatrix}1\\0\end{pmatrix} , \begin{pmatrix}1\\0\end{pmatrix}>.$$ Here we want to use the fact that $\phi$ has a certain amount of regularity therefore the quantity  $E\Delta P$ has small area. 

$$|<I_{22} \begin{pmatrix}1\\0\end{pmatrix} , \begin{pmatrix}1\\0\end{pmatrix}>|= |\frac{\omega_0}{2\pi} \int_{E\Delta P \cap \{ d(x_0)\leq |y| \leq \delta \}} \frac{y_{1}y_{2}}{|y|^2} G(|y|) dy|$$ $$ \leq \frac{\omega_0}{2\pi} \int_{E\Delta P \cap \{ |y-(0,d(x_0))| \leq 2\delta \}}  G(|y|) dy,
$$

since $\delta \leq d(x_0).$

Recall that $\partial E$ is of class $C^{1+\mu}.$ Then, we can expect to bound the measure of $E\Delta P$ using the regularity of $\partial E.$

Given $\epsilon>0,$ it is easy to see that since $P$ is the tangent plane to $E$ at the point $(0,d(x_0)),$ we have:

$$|E\Delta P \cap \{|y-(0,d(x_0))| \leq \epsilon\}| \lesssim  \epsilon^{2+\mu} |\phi|_{C^\mu}.$$

Therefore,

$$|<I_{22} \begin{pmatrix}1\\0\end{pmatrix} , \begin{pmatrix}1\\0\end{pmatrix}>|\lesssim |\phi|_{\mu} \int_{0}^{2\delta} \epsilon^{1+\mu}G(\epsilon)d\epsilon \leq |\phi|_{\mu}\int_{0}^{2\delta} \epsilon^{\mu-1} \text{Log}(1+\frac{1}{\epsilon})d\epsilon. $$

Now make a change of variables and we see:

$$|<I_{22} \begin{pmatrix}1\\0\end{pmatrix} , \begin{pmatrix}1\\0\end{pmatrix}>|\lesssim |\phi|_{\mu} \delta^\mu\int_{0}^{2} \gamma^{\mu-1} \text{Log}(1+\frac{1}{\delta\gamma})d\gamma \leq C(\mu)Log(1+\frac{1}{\delta}). $$

This tells us that $<I_{2} \begin{pmatrix}1\\0\end{pmatrix} , \begin{pmatrix}1\\0\end{pmatrix}>$ is bounded by $\text{Log}(1+\Delta_\mu)$ up to the error term which we deal with next. 

\subsection{Dealing with the Error Term}
We now need to consider the term $$|I_{2} \frac{m(\Delta_\mu)}{m(\frac{1}{d(x_0)})}|.$$

To estimate this term we employ the original method of Constantin-Bertozzi.

For every $\rho\geq d(x_0)$, consider the directions $z$ so that $x_0 +\rho z \in E:$

$$S_{\rho}(x_{0})= \{z :  |z|=1, x_0 +\rho z \in E \}.$$

Choose a point $\tilde{x} \in \partial E$ so that $|x_0 -\tilde{x}|=d(x_0)$ and consider the semicircle
$$\Sigma(x_0)=\{z: |z|=1, \nabla\phi(\tilde{x})\cdot z \geq 0 \}.$$ 

As $d(x_0)$ approaches zero, $R_{\rho}(x_0)=S_{\rho}(x_0) \Delta \Sigma_{\rho}(x_0)$ becomes small. The following lemma is taken from [4]:

\begin{lemm}
If $H^1$ denotes the Lebesgue measure on the unit circle then 

$$H^{1}(R_\rho (x_0)) \leq 2\pi ((1+2^\mu)\frac{d(x_0)}{\rho}+2^\mu (\frac{\rho}{\mu})^\mu)$$ for all $\rho\geq d(x_0), \mu>0, x_0$ so that $d(x_0)<\delta=(\frac{|\nabla \phi|_{\inf}}{|\nabla \phi|_\mu})^{1/\mu}.$

\end{lemm}

If we write $|I_{2} \frac{m(\Delta_\mu)}{m(\frac{1}{d(x_0)})}|$ in polar coordinates we get that it is majorized in the following way:

$$|I_{2} \frac{m(\Delta_\mu)}{m(\frac{1}{d(x_0)})}|\leq \frac{m(\Delta_\mu)}{m(\frac{1}{d(x_0)})}\frac{|\omega_0|}{2\pi}\int_{d(x_0)}^{\delta} \frac{d\rho}{\rho}m(\frac{1}{\rho})H^{1}(R_\rho (x_0)).$$

$$\leq m(\Delta_\mu)\frac{|\omega_0|}{2\pi}\int_{d(x_0)}^{\delta}\frac{d\rho}{\rho}H^{1}(R_\rho (x_0))$$

Now using the lemma on the size of $H^1(R_\rho(x_0)),$ we see that 

$$|I_{2} \frac{m(\Delta_\mu)}{m(\frac{1}{d(x_0)})}| \leq C(\mu)(1+m(\Delta_\mu)\text{Log}(\Delta_\mu)).$$

This concludes the proof of Lemma 6.1.

\qed

\section{Estmating the main Non-Linear Term}

$$\nabla^{\perp}\phi_t +(u\cdot \nabla) \nabla^{\perp} \phi= \nabla u \nabla^{\perp}\phi$$

The goal of this section is to show how we can deal with the non-linearity $\nabla u \nabla^\perp \phi.$ Without using any structural features of our problem, classical estimates yield:

\begin{equation} |\nabla u \nabla^{\perp}\phi|_\mu \leq C(|\nabla u|_{\infty} |\nabla \phi|_\mu+|\nabla u|_{\mu} |\nabla \phi|_\infty). \end{equation} Unfortunately, $\nabla u$ does not have any Holder regularity. Indeed, $\nabla u$ is not even bounded for our model. Fortunately, we can make use of a structural feature of our problem, as was noticed by Chemin and Bertozzi-Constantin, to put all of the regularity onto $\nabla \phi.$ The following estimate was proved for the classical problem: \begin{equation} |\nabla u \nabla^{\perp}\phi|_\mu \leq C|\nabla u|_{\infty} (|\nabla \phi|_\mu+1). \end{equation} Unfortunately, in our case such an estimate still won't work for two reasons:

\vspace{2mm}

(1) $\nabla u$ is not bounded in general.

\vspace{2mm}

(2) The insertion of the $m(|D|)$ term in the Biot-Savart law kills the classical identities and estimates that led to (7.2).

\vspace{2mm}

Fortunately for us, we are able to prove the following estimate using ideas along the lines of the original work of Bertozzi-Constantin along with a commutator estimate which allows us to work in basically the classical regime.
Indeed, define $$u_{classical}(x,t)=\tilde{u}(x,t)=\frac{a_0}{2\pi}\int_{E}\nabla^{\perp}_{x} (\text{Log}|x-y|)dy.$$ If $W$ is a divergence free vector field tangent to the boundary of $E$ then recall (see Proposition 2 of [4]) that \begin{equation} \nabla \tilde{u} W(x)= \int_{E} \frac{ \sigma(|x-y|)}{|x-y|^2}(W(x)-W(y))dy \end{equation}.

We can write 

$$\nabla u W= (m(|D|)\nabla \tilde{u}) W= m(|D|) (\nabla\tilde{u}W)+[m(|D|),\nabla \tilde{u}]W$$

Now, by virtue of (7.1), the following estimate holds for the first term on the right hand-side above:

$$\sup_{j} 2^{\sigma j}\frac{1}{m(2^j)}|\Delta_{j} m(|D|) (\nabla\tilde{u}W)|_{\infty} \leq |\nabla\tilde{u}W|_{C^\sigma} \leq |\nabla\tilde{u}|_{\infty}|W|_{C^\sigma}. $$

The first inequality is obvious. The second uses (7.1).

Now we must deal with the commutator $[m(|D|),\nabla \tilde{u}]W.$ The estimate for the commutator is given in Lemma 5.2.

Indeed, it is possible to show that 

$$\sup_j 2^{\sigma j} \frac{1}{m(2^j)}|\Delta_{j}([m(|D|),\nabla \tilde{u}]W)|_{\infty} \leq |\nabla \tilde{u}|_{\infty}|W|_{C^\sigma}$$

Putting all this together we get the following lemma:
\begin{lemm}
Assume $u$ is given by the Modified Biot-Savart law:

$$u(x,t)= \frac{a_{0}}{2\pi}m(|D|)\int_{E}\nabla^{\perp}_{x} (\text{Log}|x-y|)dy.$$

Where $$E = \{ x\in \mathbb{R}^2 | \, \phi(x) >0 \}.$$ Let $W$ be a divergence free vector field which is tangent to the boundary of $E$.

Then,

$$\sup_{j}2^{j\sigma}\frac{1}{ m(2^j)} |\nabla u W| \leq |\nabla \tilde{u}|_{\infty}|W|_{\sigma} \leq C(\sigma)(1+\text{Log}(\Delta_\sigma))|W|_{\sigma}.$$

\end{lemm}
\qed 

\section{Steps 5 and 6}

Step 1 tells us that $u$ is not quite Lipschitz, but it is better than quasi-Lipschitz. In fact, the bound:

$$\sup_{j} \frac{S_j \nabla u}{m(2^j)} \leq C$$ is equivalent to the bound: $$\sup_{x\neq y} \frac{|u(x)-u(y)|}{|x-y|(1+m(\frac{1}{|x-y|}))} < C$$ (see [2]). 

 This tells us that even though we cannot propogate the $C^\mu$ regularity of $\nabla \phi$ in transporting by $u,$ we can show that we only lose an arbitrarily small amount of regularity. 
 
Indeed, since $\nabla u=m(|D|)\nabla \tilde{u},$ and since (from the classical estimates) $|\nabla \tilde{u}|_\infty \leq C(\mu)\text{Log}(1+\Delta_\mu)$ we easily see  

$$|S_j \nabla u|_\infty \leq C m(2^j)|\nabla \tilde{u}|_\infty \leq C m(2^j) (1+\text{Log}(\Delta_\mu)).$$

Apply $\Delta_j$ to (4.2). Then we get:

$$ \Delta_j (\nabla^{\perp} \phi )_t +(u\cdot \nabla) \Delta_j (\nabla^{\perp} \phi )= \Delta_j(\nabla u \nabla^{\perp} \phi) +CO_j$$

Where $CO_j$ is the following commutator:

$$CO_j = [u\cdot \nabla, \Delta_j] \nabla^{\perp}\phi.$$

Using the bound (4) above, $$  \frac{1}{m(2^j)}2^{\sigma j}|CO_j|_{\infty} \leq  (C(\mu)(1+\text{Log}(\Delta_{\mu})) | (\nabla^{\perp}\phi)|_{\sigma},$$ with bounds independent of $j.$

\section{The Losing Estimates}

Take (4.2) and dot it with $\nabla^\perp \phi.$ 
Using the transport structure we easily get:

\begin{equation} (|\nabla^\perp \phi|^2)_{t}+ (u\cdot \nabla) |\nabla^\perp \phi|^{2}= <\nabla u\nabla^\perp \phi, \nabla^\perp \phi>=<\nabla u \tau,\tau>|\nabla^\perp\phi|^2,\end{equation} where $\tau=\frac{\nabla^\phi}{|\nabla^\phi|}$ as in section 6. Lemma 6.1 tells us that $<\nabla u \tau,\tau>$ is bounded which is the crucial. 

This implies that

$$ \text{Log}(|\nabla^\perp \phi|\circ \Phi)_t \leq  |<\nabla u\tau,\tau>\circ \Phi|,$$ where $\Phi$ is the flow-map associated to $u.$

Thus we get:

$$|\nabla \phi|_{\inf} (t) \geq |\nabla \phi_0|_{\inf} \, \text{exp}({-\int_{0}^{t}|<\nabla u \tau,\tau>|_{\infty}(s)ds}). $$

So that

\begin{equation} |\nabla \phi|_{\inf} (t) \geq |\nabla \phi_0|_{\inf} \, \text{exp}({-\int_{0}^{t}C(\mu)m(\Delta_\mu)(1+\text{Log}(\Delta_\mu))(s)ds}).  \end{equation}

Note that (9.2) is true for any $\mu>0.$

We now need to show a similar bound on the evolution of the $|\nabla \phi|_\mu .$ We will then be able to conclude using Osgood's lemma.

Fix $T>0.$  Define $V(t)=\int_{0}^{t}C(\gamma)(1+\text{Log}(\Delta_\gamma))(s)ds$ for some $\gamma>0$ fixed ($\gamma$ will be taken to be $\mu-\epsilon$). Take $\eta=\frac{\epsilon}{V(T)}$.

Let $\mu_t = \mu -\eta V(t).$

$$ \Delta_j (\nabla^{\perp} \phi )_t +(u\cdot \nabla) \Delta_j (\nabla^{\perp} \phi )= \Delta_j(\nabla u \nabla^{\perp} \phi) +CO_j$$ 

This implies that $$|\Delta_j \nabla \phi|_\infty \leq |\Delta_j \nabla \phi_0|_\infty+ \int_{0}^{t} |\Delta_j(\nabla u \nabla^{\perp} \phi)|_\infty +|CO_j|_\infty$$

so that $$2^{(j+2)\mu_t }|\Delta_j \nabla \phi|_\infty \leq 2^{(j+2)\mu}|\Delta_j \nabla \phi_0|_\infty+ \int_{0}^{t} 2^{(j+2)\mu_t }(|\Delta_j(\nabla u \nabla^{\perp} \phi)|_\infty (s) +|CO_j|_\infty (s) )ds,$$
 $$2^{\mu_t j}|\Delta_j \nabla \phi|_\infty \leq 4|\nabla \phi_0|_\mu \int_{0}^{t} 2^{(j+2)(\mu_t-\mu_s)} 2^{(j+2)\mu_s }(|\Delta_j(\nabla u \nabla^{\perp} \phi)|_\infty (s) +|CO_j|_\infty (s) )ds,$$
and $$2^{(j+2)\mu_t}|\Delta_j \nabla \phi|_\infty \leq 4|\nabla \phi_0|_\mu+ \int_{0}^{t} 2^{(j+2)(\mu_t-\mu_s)}m(2^j) 2^{(j+2)\mu_s }\frac{1}{m(2^j)}(|\Delta_j(\nabla u \nabla^{\perp} \phi)|_\infty (s) +|CO_j|_\infty (s) )ds.$$

Now using the estimates on $\nabla u W,$ $W=\nabla^\perp \phi$ and $CO_j,$ we get:

$$2^{(j+2)\mu_t}|\Delta_j \nabla \phi|_\infty \leq 4|\nabla \phi_0|_\mu+ \int_{0}^{t} 2^{(j+2)(\mu_t-\mu_s)}m(2^j) C(\gamma)(1+\text{Log}(\Delta_\gamma))|\nabla\phi|_{\mu_s}ds,$$
$$=|\nabla \phi_0|_\mu+ \int_{0}^{t} 2^{-\eta(j+2)(\int_{s}^{t}C(\gamma)(1+\text{Log}(\Delta_\gamma))(s)ds)}m(2^j) C(\gamma)(1+\text{Log}(\Delta_\gamma))|\nabla\phi|_{\mu_s}ds  $$
for any $\gamma>0.$ Now we want to break our analysis into two cases. The first case being when $j$ is large in a certain sense and the second when $j$ is small. 

\emph{Case 1:} $$ \frac{j+2}{m(2^j)} \geq \frac{C}{\eta}.$$

If $C$ is chosen large enough, we see that $$\int_{0}^{t} 2^{-\eta(j+2)(\int_{s}^{t}C(\gamma)(1+\text{Log}(\Delta_\gamma))(\alpha)d\alpha)}m(2^j) C(\gamma)(1+\text{Log}(\Delta_\gamma))ds \leq \frac{1}{2}$$ (note that the integral $ds$ can be computed exactly).
Thus we see that $$ 2^{(j+2)\mu_t}|\Delta_j \nabla \phi|_\infty \leq 4|\nabla \phi_{0}|_{\mu} +\frac{1}{2}\sup_{t} |\nabla \phi|_{\mu_{t}}.$$

Once more, this is only true for $j$ large. 

\emph{Case 2:} $$ \frac{j+2}{m(2^j)} \leq \frac{C}{\eta}.$$

For these values of $j,$ (keeping the example of $m=\text{log log}$ in mind, $$m(2^j) \leq m(2^\frac{C}{\eta}) \leq Cm(2^\frac{1}{\eta}).$$

Now, $$2^{(j+2)\mu_t}|\Delta_j \nabla \phi|_\infty \leq 4|\nabla \phi_0|_\mu+ \int_{0}^{t} 2^{(j+2)(\mu_t-\mu_s)}m(2^j) C(\gamma)(1+\text{Log}(\Delta_\gamma))|\nabla\phi|_{\mu_s}ds.$$

$$ \leq 4|\nabla \phi_0|_\mu+ \int_{0}^{t} m(2^\frac{1}{\eta}) C(\gamma)(1+\text{Log}(\Delta_\gamma))|\nabla\phi|_{\mu_s}ds$$
$$=4|\nabla \phi_0|_\mu+ \int_{0}^{t} m(2^\frac{V(T)}{\epsilon}) C(\gamma)(1+\text{Log}(\Delta_\gamma))|\nabla\phi|_{\mu_s}ds.$$

So, for $j$ large we have:

$$ 2^{(j+2)\mu_t}|\Delta_j \nabla \phi|_\infty \leq 4|\nabla \phi_{0}|_{\mu} +\frac{1}{2}\sup_{t} |\nabla \phi|_{\mu_{t}}.$$

And for $j$ small we have:

$$2^{(j+2)\mu_t}|\Delta_j \nabla \phi|_\infty \leq 4|\nabla \phi_0|_\mu+ \int_{0}^{t} \frac{1}{\epsilon}m(2^{V(T)}) C(\gamma)(1+\text{Log}(\Delta_\gamma))|\nabla\phi|_{\mu_s}ds.$$

Putting these two together we get:

$$|\nabla \phi|_{\mu_{T}} \leq 2|\nabla \phi_{0}|_\mu + \frac{1}{\epsilon}m(2^{V(T)})C(\gamma)\int_{0}^{T}(1+\text{Log}(\Delta_\gamma))|\nabla\phi|_{\mu_s}ds.$$

Now apply the Gronwall lemma and we see:

\begin{equation}|\nabla \phi|_{\mu-\epsilon}(T) \leq 2|\nabla \phi_{0}|_\mu \text{exp}[\frac{C(\gamma)}{\epsilon}m(\text{exp(V(T))})V(T)],\end{equation}

where, again, $V(t)=\int_{0}^{t}C(\gamma)(1+\text{Log}(\Delta_\gamma))(s)ds.$

Now recall equation (9.2):

$$|\nabla \phi|_{\inf} (t) \geq |\nabla \phi_0|_{\inf} \, \text{exp}({-\int_{0}^{t}C(\mu)m(\Delta_\mu)(1+\text{Log}(\Delta_\mu))(s)ds}).$$

Take $\gamma=\mu-\epsilon.$ 

Recall that $\Delta_\gamma=\frac{|\nabla \phi|_\gamma}{|\nabla \phi|_{\inf}}.$ 

So we see that the quantity $\Delta:=\Delta_{\mu-\epsilon}$ satisfies the following inequality:

$$\text{Log}(\Delta)(T) \leq C_{0}+C[m(\text{exp}(\int_{0}^{T}\text{Log}(\Delta)))\int_{0}^{T}\text{Log}(\Delta)+\int_{0}^{T}m(\Delta)(1+\text{Log}(\Delta))].$$ Where $C_{0}$ is a constant depending only on $ |\phi_{0}|_\mu, |\phi_{0}|_{\inf}$ and C depends only upon $\epsilon$ and $\mu.$

Let $\tilde{m}(\cdot)=m(\text{exp}(\cdot)).$ 
So that $f:=\text{Log}(\Delta)$ satisfies:

$$f(T) \leq C_{0}+C[\tilde{m}(\int_{0}^{T}f)\int_{0}^{T}f+\int_{0}^{T}\tilde{m}(f)(1+f)].$$

Since $\tilde{m}(\cdot)$ grows like $\text{Log}(1+\cdot)$ (or just a little more), we see that we can apply lemma 1.2 to conclude that $f$ must remain bounded for all time.
This concludes the proof of Theorem 3.1.
\qed

\section{Acknolwedgements}

The author acknowledges the support of NSF grant DMS-0807347. The author wishes to thank his advisor N. Masmoudi for his continued support throughout the completion of this project. We also thank professor A. Kiselev for many helpful comments. Finally, the author thanks the anonymous referee for his/her very helpful comments.

\end{document}